\documentclass[reqno,11pt]{amsart}
\usepackage{graphicx}
\usepackage{mathrsfs}
\usepackage{a4}
\usepackage{amsmath}
\usepackage{amsfonts}
\usepackage{amssymb}

\newtheorem{thm}{Theorem}[section]

\newtheorem{prop}[thm]{Proposition}
\theoremstyle{definition}
\newtheorem{defn}[thm]{Definition}
\theoremstyle{remark}
\newtheorem{rem}[thm]{Remark}
\numberwithin{equation}{section}

\begin{document}
\title[Stability of JS-queue models]
{The stability of join-the-shortest-queue models with general input and output processes}
\author[ABRAMOV]{Vyacheslav M. Abramov}
\address{School of Mathematical Sciences, Monash University, Clayton campus, Building 28,
Level 4, Wellington road, Clayton, Victoria 3800, Australia}
\email{Vyacheslav.Abramov@sci.monash.edu.au} 

\begin{abstract}
The paper establishes necessary and sufficient conditions for
stability of different join-the-shortest-queue models including
load-balanced networks with general input and output processes. It
is shown, that the necessary and sufficient condition for
stability of load-balanced networks is associated with the
solution of a linear programming problem precisely formulated in
the paper. It is proved, that if the minimum of the objective
function of that linear programming problem is less than 1, then
the associated load-balanced network is stable.
\end{abstract}

\maketitle

\newpage
%
%
%
%

%
%
%
%

\section{Introduction}

\subsection{The goal of the paper}
The theory of parallel queues is a distinguished area of queueing
theory. Parallel queues have good properties (e.g. Halfin
\cite{Halfin (1985)}, Mitzenmacher \cite{Mitzenmacher -
dissertation}, Winston \cite{Winston (1977)}) resulting in various
applications in different areas of science and technology. The
literature on parallel queues is very rich and includes solutions
of a large number of various theoretical and applied problems. In
this context we mention the pioneering paper of Flatto and McKean
\cite{Flatto and McKean (1977)}, the papers of Flatto \cite{Flatto
(1985)}, Flatto and Hahn \cite{Flatto and Hahn (1984)}, as well as
the recent paper of Kurkova and Suhov \cite{Kurkova Suhov (2003)}
for mentioning a few.

During the last years there has been increasing interest in
parallel queueing systems due to a growing development of
telecommunication technologies. As a result, there has been
substantially increasing number of related publications including
such areas as routing and admission control problems, scheduling
and many other areas using the join-the-shortest-queue (JS-queue)
policy.

Despite of many existing publications on JS-queue models, it is
still not a well-studied area. There is a large number of unsolved
problems or problems the solution of which is very far from its
completion. One of them is a problem of stability. The stability
(instability) of different JS-queue (together with stability other
closely related models) has been studied in many papers. We refer
to the papers of Foley and McDonald \cite{Foley McDonald (2001)},
Foss and Chernova \cite{Foss and Chernova (1991)}, \cite{Foss and
Chernova (1998)}, Kurkova \cite{Kurkova (2001)}, Sharifnia
\cite{Sharifnia (1997)}, Suhov and Vvedenskaya \cite{Suhov and
Vvedenskaya (2002)}, Tandra, Hemachandra and Manjunath
\cite{Tandra Hemachandra and Manjunath (2004)}, Vvedenskaya and
Suhov \cite{Vvedenskaya and Suhov (2004)},
 and Vvedenskaya,
Dobrushin and Karpelevich \cite{Vvedenskaya Dobrushin and
Karpelevich (1996)}. Particularly, the stability of a Markovian
load-balanced network, consisting of two stations only, has been
studied in \cite{Kurkova (2001)}.





In the present paper, we establish necessary and sufficient
conditions of stability for a wide class of load-balanced
networks, the arrival and departure processes of which are assumed
to be quite general. The exact description of these processes will
be provided later in the paper.

\subsection{Motivation} The stability of stochastic processes and
especially queueing systems and networks of queues is a very
significant area of research. The first results on the stability
of queues go back to the well-known pioneering papers of Lindley
\cite{Lindley (1952)} and Kiefer and Wolfowitz \cite{Kiefer and
Wolfowitz (1955)} on the stability of classic single-server and
multiserver queueing systems. One of the most significant
contribution to stability of queueing systems was made slightly
later by Loynes \cite{Loynes (1962)}, who was the first person to
establish the stability of single-server queues with dependent
interarrival and service times. The approach of Loynes
\cite{Loynes (1962)} is well-known and has been widely used to
prove the stability of various queueing systems. It led to a new
vision of the stability for models with \textit{dependent}
interarrival and/or service times, and has been a source for new
methods of the stability of complex queueing models. Among them
there are renovative theory and recurrence equation methods (e.g.
\cite{Akhmarov and Leont'eva (1976)}, \cite{Altman and Borovkov
(1997)}, \cite{Altman and Hordijk (1997)}, \cite{Asmussen and Foss
(1993)}, \cite{Borovkov (1984)}, \cite{Borovkov (1998)},
\cite{Brandt Franken and Lizek (1990)}, \cite{Foss and Kalashnikov
(1991)}
 and many others) as well as a martingale-based approach of \cite{Abramov
 (2006)}.

 However, the development of the Loynes approach towards queueing
 networks, because of their complicated nature, has been problematical, and the proof of the stability and
 ergodicity of queueing networks
is typically based on other methods based on the special theories
of Markov and regeneration processes.

The first results on the stability of Jackson-type queueing
networks have been obtained by Borovkov \cite{Borovkov
 (1986)}, \cite{Borovkov (1987)}. Then the stability and
 ergodicity of networks have been studied in many
 papers. We mention the papers related to the stability of Jackson-type networks by Meyn and Down \cite{Meyn and Down
 (1994)}, Kaspi and Mandelbaum \cite{Kaspi and Mandelbaum (1992)},
 \cite{Kaspi and Mandelbaum (1994)}, Sigman \cite{Sigman (1990)}
and Baccelli and Foss \cite{Baccelli and Foss (1994)}. All of
these papers, but
 \cite{Baccelli and Foss (1994)}, are based on regeneration phenomena of
 the theory of Harris recurrent Markov chains.

The theory of Harris recurrent Markov chains had been exposed in
the books of Orey \cite{Orey (1971)} and Revuz \cite{Revuz
(1975)}. The detailed study of stochastic stability of Markov
chains and the theory of Harris recurrent Markov processes can be
found in the book of Meyn and Tweedie \cite{Meyn and Tweedie
 (1993)}, and in a number of research papers of these authors \cite{Meyn and Tweedie (1992-1)},
 \cite{Meyn and Tweedie (1993-2)} and \cite{Meyn and Tweedie
 (1993-3)}.


However, the proof of networks stability by the means of Harris
recurrent Markov chains is restrictive. Being very difficult
technically, it works in the cases where the sequences of
interarrival and service times consist of independent random
variables. In this case the phase space of the process can be
expanded to Markov, and a network stability is proved in terms of
the stability of the corresponding Markov process. As a rule the
proof of the network stability in this case requires specific
additional conditions. For example, in most of papers an infinite
support of interarrival time distributions is required. Sometimes
some additional method based technical conditions are required as
well (e.g. see  Dai \cite{Dai (1995)}).

Baccelli and Foss \cite{Baccelli and Foss (1994)} proved the
stability of Jackson-type networks with dependent interarrival and
service times. Their proof is based on development of renovation
theory. However, the mentioned paper is about 70 pages long, it
contain many notions and results from different areas of research
(stochastic Petri nets for example), and it is not easy to read
this paper.

In the present paper, we establish conditions for stability of
JS-queue models including load balanced networks. The class of
load-balanced networks is wider that the class of Jackson-type
networks, so the stability results of the present paper are more
general than those for Jackson-type networks. On the other hand,
our method is a Loynes-based method, and our stability results are
established for quite general networks with sequences of
\textit{dependent} interarrival and service times, and these
sequences can be dependent of one another as well. Furthermore,
for the known cases of queues and networks our results are
obtained under weaker conditions than known results. For example,
the paper of Baccelli and Foss \cite{Baccelli and Foss (1994)}
requires stationarity of the appropriate sequences of interarrival
and service times. In our case, we requires weaker conditions of
the strong law of large numbers, and our class of systems is
therefore wider.

Our challenge is as follows. We first establish the following
equivalence: \textit{the stability of usual queueing systems
follows from the stability of queueing systems with autonomous
service mechanism} (which are sometimes called queueing systems
with a walking type service \cite{Gelenbe and Iasnogorodski
(1980)}). This result is based on sample path analysis and
stochastic comparison, and although its proof is elementary, the
result is a significant contribution to the proof of the stability
of different queueing networks. Sample-path analysis for stability
of queueing systems and networks is not new (e.g. \cite{Taha and
Stidham (1999)}). However, the approach of the present paper
specifically uses a sample-path analysis in combination with other
methods and based on the new idea of reduction of the original
problem to another not traditional simpler problem.

The aforementioned stochastic inequalities established for
queueing systems is easily extended to each node of a network,
where a stable behavior of each node of a usual network is a
consequence of the stable behavior of the corresponding node of
the queueing network with autonomous service mechanism. Then the
general problem reduces to the easier problem of the stability of
networks with autonomous service, and the second part of the proof
is \textit{to establish the stability of queueing networks with an
autonomous service mechanism}. This second part of the proof is
based on the Loynes method.

Queueing systems with an autonomous service mechanism have been
introduced and initially studied by Borovkov \cite{Borovkov
(1976)}, \cite{Borovkov (1984)}, and then have been an object of
study in a large number of papers (e.g. \cite{Abramov (2000)},
\cite{Abramov (2004)}, \cite{Fricker (1986)}, \cite{Fricker
(1987)}, \cite{Gelenbe and Iasnogorodski (1980)}). In traditional
applications, queueing systems with autonomous service are
associated with a shuttle bus picking up passengers from stations.
Other, more interesting applications of these systems, are known
from computer technologies. One of such examples is the event
processing in the operating system Microsoft Windows. The details
of this issue can be found in \cite{Richter (1997)} in Section 10:
Threads Synchronization and specifically on page 396 (Events). The
original construction there is much more complicated than that in
our example described below, and it is explained in terms of
threads and their synchronization which is required in the Windows
programming. However, loosely speaking, it can be explained as
follows. In specified time instants, the event processor checks
whether there is an event (such as mouse-movement, mouse-click,
mouse-double-click and so on) in the event queue. If such an event
exists (there is a thread receiving an acknowledgement (signal)
about it), then the system processes it at specified time instant.
Otherwise, the system continues to check the state of the event
queue.

\smallskip

The idea of reducing one stability problem of a complex network to
another corresponding stability problem of a network with
simpler/concise properties is not new. There are special criteria
in the literature allowing to replace an original (stochastic)
network by its (deterministic) fluid model to study the stability.
Such criteria for quite general class of queueing networks with
multiple customer classes has been established by Dai \cite{Dai
(1995)}. Formally it had been used before for establishing the
instability of a special configuration of a network with priority
classes by Rybko and Stolyar \cite{Rybko and Stolyar (1992)}.

Although the reducing an original stochastic network to its
deterministic fluid analogue looks natural, in fact it requires
additional (mild) conditions. In the case of our study by reducing
an original stochastic network to its associated stochastic
network with autonomous service mechanism no additional condition
is required.

\subsection{Organization of the paper} The rest of the paper is
organized as follows. In Section 2 we describe the main JS-queue
models, which will be then developed and modified in the following
sections. (The material of the paper is presented in the order of
increasing complexity.) In the same section we give all of the
necessary definitions related to stability of queues and networks.
In Section 3 we prove the correspondence between the stability of
the original queueing system and that of the associated queueing
system with an autonomous service mechanism. The proof is based on
sample path analysis. In Section 4 we establish conditions for
stability of JS-queue models of queueing systems, and then in
Section 5 we establish conditions for stability of load-balanced
networks. In Section 6 we conclude the paper, where the stability
of more general networks, than those studied here, with batch
arrival and service times are discussed.

\section{Description of the main models and definitions of stability}
\subsection{Main models}
In this section we describe main JS-queue models with an
autonomous service mechanism. These models and some of the
assumptions related to the arrival and departure processes will be
then modified in the following sections.

\smallskip
$\bullet$ There are $m$ identical servers, each of which having
its own queue.

\smallskip
$\bullet$ All of the processes that describe queueing models are
assumed to be right-continuous and to have left limits.

\smallskip
$\bullet$ The arrival process is governed by two point processes
$A(t)$ and $A^{\prime}(t)$. The process $A(t)$ is defined by
sequence $\{\tau_{n}\}_{n\geq 1}$ of positive random variables,
and the corresponding sequence of points is the following:
$t_{1}=\tau_{1}$, and $t_{n+1}=t_{n}+\tau_{n+1}$, $n\geq 1$. Then,
$A(t)=\sum_{i=1}^{\infty }\mathbf{1}_{\{t_{i}\leq t\}}$. The
process $A^{\prime}(t)$ is defined analogously. We have the
sequence of positive random variables
$\{\tau_{n}^{\prime}\}_{n\geq 1}$,
and the sequence of points
$t_{1}^{\prime}=\tau_{1}^{\prime}$, and $t_{n+1}^{\prime}=t_{n}^{\prime}%
+\tau_{n+1}^{\prime}, n\geq 1$. Then,
$A^{\prime}(t)=\sum_{i=1}^{\infty
}\mathbf{1}_{\{t_{i}^{\prime}\leq t\}}$.
We assume
\begin{equation}\label{r1.1}
\mathbf{P}\Big\{\lim_{t\rightarrow\infty}\frac{A(t)}{t}=\lambda\Big\}%
=1,
\end{equation}
and
\begin{equation}\label{r1.2}
\mathbf{P}\Big\{\lim_{t\rightarrow\infty}\frac{A^{\prime}(t)}{t}%
=\lambda^{\prime}\Big\}=1.
\end{equation}
The process $A(t)$ forms a \textit{dedicated} traffic, while the
process $A^{\prime}(t)$ forms an \textit{opportunistic} traffic.

\smallskip
$\bullet$ A customer arriving at moment $t_{n}$, $n\ge1$, is
assigned to the
$j$th queue, $j=1,2,\ldots,m$, with the probability $p_{j}$ ($\sum_{j=1}^{m}%
p_{j}=1$), residing there to wait for the service.

\smallskip

$\bullet$ A customer, arriving at moment $t_{n}^{\prime}$,
$n\ge1$, is assigned to the queue with the shortest queue-length
breaking ties at random.

\smallskip
$\bullet$ The departure process from the $j$th server is governed
by the point process $D^{(j)}(t)$. Specifically, the $n$th service
time of the $j$th server is denoted $\chi_{n}^{(j)}$, and the
corresponding sequence of points is denoted $\{x_n^{(j)}\}$ where
$x_1^{(j)}=\chi_1^{(j)}$ and
$x_{n+1}^{(j)}=x_{n}^{(j)}+\chi_{n+1}^{(j)}$, $n\geq1$. We assume
\begin{equation}\label{r1.3}
\mathbf{P}\Big\{\lim_{t\rightarrow\infty}\frac{D^{(j)}(t)}{t}=\mu\Big
\}=1.
\end{equation}
%
%
%

$\bullet$ For our convenience we assume that the processes $A(t)$,
$A^{\prime }(t)$ and $D^{(j)}(t)$, $j=1,2,\ldots,m$, all are
mutually independent point processes. (Then this condition
together with other conditions \eqref{r1.1}, \eqref{r1.2} and
\eqref{r1.3} will be relaxed.)

\smallskip

$\bullet$ The service mechanism of each server is assumed to be
autonomous. This means the following. Let $Q^{(j)}(t)$ denote the
number of customers in the $j$th queue at time $t$,
$j=1,2,\ldots,m$, and let $Q^{(j)}(0)=0$. Let $A^{(j)}(t)$ and
$A^{\prime}{}^{(j)}(t)$ denote the thinning of the processes
$A(t)$ and $A^{\prime}(t)$ respectively, where $A^{(j)}(t)$ and
$A^{\prime}{}^{(j)}(t)$ are arrival processes to the $j$th queue.
Then,
\begin{equation}\label{r1.4}
Q^{(j)}(t)=A^{(j)}(t)+A^{\prime}{}^{(j)}(t)-\int_{0}^{t}\mathbf{1}%
_{\{Q^{(j)}(s-)>0\}}\mbox{d}D^{(j)}(s).
\end{equation}
For the further convenience the above model is denoted $\wp_{m}$,
where the subscript $m$ denotes the number of parallel queues.

\smallskip
The model $\wp_{m}$ is a special case of the more general model,
in which it is assumed that there are $m$ different arrival point
processes $A^{(1)}(t)$, $A^{(2)}(t)$,\ldots, $A^{(m)}(t)$ of
dedicated traffic corresponding to $m$ servers. Let us denote this
more general model $\Gamma_m$.

\smallskip
In turn, we will also consider the particular case of the model
$\wp_{m}$, where
$$p_{1} = p_{2} = \ldots = p_{m} = \frac{1}{m}.$$
In this
case the families $\{A^{(j)}(t)\}_{j\leq m}$ and $\{A^{\prime
(j)}(t)\}_{j\leq m}$ consist of identically distributed processes.
The above symmetric model with $m$ parallel queues is denoted
$\Sigma_{m}$.

\subsection{Definitions of stability}
For the sake of convenience we discuss definitions for $\wp_{m}$
models. The extension of this definition to models $\Gamma_m$ is
technical. The above equation for model $\wp_{m}$ is given for all
$t>0$. For our purpose we extend this equation, assuming that all
the processes start at $a$. Then, instead of \eqref{r1.4} we have
the following equation:
\begin{equation}\label{r1.5}
Q_a^{(j)}(t)=A_a^{(j)}(t)+A_a^{\prime}{}^{(j)}(t)-\int_{a}^{t}\mathbf{1}%
_{\{Q_a^{(j)}(s-)>0\}}\mbox{d}D_a^{(j)}(s),
\end{equation}
where the subscript $a$ says that the processes start at $a$.
\begin{defn}\label{defn2}
The system $\wp_{m}$ is called to be instable if
$$\lim_{a\to~-\infty}\mathbf{P}\{Q_{a}^{(j_0)}(t)\in\mathscr{S}\}=0$$
for any bounded set $\mathscr{S}$ at least for some index $j_0$
and any $t$. Otherwise, the system is called to be stable.
\end{defn}

Then, the stability of the system means the following

\begin{defn}\label{defn1}
The system $\wp_{m}$ is said to be stable if there exists a
bounded set $\mathscr{S}$ such that
$$\limsup_{a\to~ -\infty}\mathbf{P}\{Q_{a}^{(j)}(t)\in\mathscr{S}\}>0$$
for all $j=1,2,\ldots,m$ and any $t$.
\end{defn}

This definition remains in force for all of JS-queue models
included load-balanced networks considered in the paper.

\smallskip
In some examples, the processes all are assumed to start at zero.
In this case by the stability of the system $\wp_{m}$ we mean the
existence of a bounded set $\mathscr{S}$ such that
$$\limsup_{t\to\infty}\mathbf{P}\{Q^{(j)}(t)\in\mathscr{S}\}>0$$
for all $j=1,2,\ldots,m$.

%
%

\section{Sample-path comparison of queueing systems}

In this section we compare three different queueing systems given
on the same probability space $(\Omega, \mathscr{F}, \mathbf{P})$
and therefore defined by the same governing sequences of random
variables, but different specific rules of departure. For the sake
of simplicity we assume that all of these three systems start at
zero with empty queue.

These systems are defined by an arrival point process $A(t)$ and
departure process $D(t)$. These processes are defined by the
corresponding governing sequences $\{\tau_n\}$ and $\{\chi_n\}$.
Let $t_k=\tau_1+\tau_2+\ldots+\tau_k$ and
$x_k=\chi_1+\chi_2+\ldots+\chi_k$. Then, the point processes
$A(t)$ and $D(t)$ are
$$
A(t)=\sum_{i=1}^\infty\mathbf{1}_{\left\{t_i\leq t\right\}}, \ \
D(t)=\sum_{i=1}^\infty\mathbf{1}_{\left\{x_i\leq t\right\}}.
$$
The first queueing system is the queueing system with autonomous
service, which is denoted $\mathscr{Q}_1$. The queue-length
process for this system is defined as
$$
Q_1(t)=A(t)-\int_0^t\mathbf{1}_{\{Q_1(s-)>0\}}\mbox{d}D(s).
$$

The second queueing system is the usual queueing system. Denote
this system $\mathscr{Q}_2$, and the queue-length process of this
queueing system $Q_2(t)$ is defined by the traditional recurrence
equations well-known from the queueing theory. Specifically,
following \cite{Borovkov (1976)}, p.19, the queue-length process
$Q_2(t)$ is defined by \textit{interrupted governing sequences} as
follows. Denoting $\eta_1=\inf\{k: t_1+x_k<t_{k+1}\}$ we have the
following relations for the queue-length process $Q_2(t)$. For
$0\leq t<t_1$, we have $Q_2(t)=0$, and for $t_1\leq t<
t_1+x_{\eta_1}$, the queue-length $Q_2(t)$ is the difference
between the number of arrivals and service completions during the
interval [$t_1, t$] including the arrival at the instant $t_1$.
Then, for $t_1+x_{\eta_1}\leq t< t_{\eta_1+1}$ we have $Q_2(t)=0$.
Next, denoting $\eta_2=\inf\{k>\eta_1:
t_{\eta_1+1}+x_k-x_{\eta_1}<t_{k+1}\}$, then for $t_{\eta_1+1}\leq
t<t_{\eta_2+1}$ the queue-length $Q_2(t)$ is the difference
between the number of arrivals and service completions during the
interval [$t_{\eta_1+1},t$) including the arrival at the instant
$t_{\eta_1+1}$. The stopping times $\eta_3$, $\eta_4$, \ldots are
defined similarly.

The third queueing system is a special queueing system with
delayed departures is denoted $\mathscr{Q}_3$. The queue-length
process of this system is defined as follows. The arrival process
$A(t)$ is the same as in the systems $\mathscr{Q}_1$ (or
$\mathscr{Q}_2$), but departures of the customers are delayed as
follows. The instant of first departure is $\chi_1+\chi_2$, the
instant of the second departure is $\chi_1+\chi_2+\chi_3$ and so
on. The departures occur only if there is at least one customer in
the system. Thus, the difference between queueing systems
$\mathscr{Q}_1$ and $\mathscr{Q}_3$ is only that the departures of
customers in $\mathscr{Q}_1$ occur immediately at the above
specified time instants, while in the queueing system
$\mathscr{Q}_3$ these specified time instants are service begins
until the next specified instants correspondingly (provided that
the system $\mathscr{Q}_3$ is not empty) . The duration of these
service times are $\chi_1$ for the first customer, $\chi_2$ for
the second one, and so on. The queue-length process $Q_3(t)$ is
defined as follows. Let $\omega_A(t)$ denote the event, that the
last arrival before time $t$ was to an empty system and before
time $t$ (excluding the time instant $t$ itself) its service has
not yet been started, and $\mathbf{1}_{\omega_A(t)}$ denotes the
indicator of the event $\omega_A(t)$. Then,
\begin{equation}\label{SP1}
Q_3(t)=A(t)-\int_0^t(1-\mathbf{1}_{\omega_A(s)})\mathbf{1}_{\{Q_3(s-)>0\}}\mbox{d}D(s).
\end{equation}

According to \eqref{SP1}, the process $Q_3(t)$ is always
incremented at the moments of arrivals $t_i$, $i\geq1$. However,
it is decremented as follows. Let $t_i$ be such a moment of
arrival to an empty system (that is $Q_3(t_i)=1$), and let $\ell =
\min\{j: t_i\leq x_j\}$. Then at the moment $x_\ell$ the
queue-length is not decremented, i.e. $Q_3(x_\ell)=Q_3(t_i)$. In
all of other points $x_j$ where the queue-length is positive and
do not satisfy the above property, the queue-length is
decremented.

\smallskip
Since all of the queue-length processes are defined on the same
probability space, then these processes are provided by the
additional argument $\omega\in\Omega$ in the places where it is
required.
\begin{prop}\label{prop1}
\begin{equation}\label{3.1}
Q_2(t,\omega)\leq Q_3(t,\omega).
\end{equation}
\end{prop}
\begin{proof}
Using sample path analysis we prove this proposition as follows.
In the system $\mathscr{Q}_3$ we do not consider the events
$\{Q_3(x_j-)=0\}$ and $\omega_A(t)$, taking into account service
completions during busy periods only. According to this note, the
points $\{x_j\}$ are to be renumbered such that there is the
correspondence of such the points between the systems
$\mathscr{Q}_3$ and $\mathscr{Q}_2$ in busy periods for the
purpose of further sample path analysis. Then it is easily seen
that the statement of the proposition follows, because in
$\mathscr{Q}_2$ the first service begin in a busy period coincides
with the moment of arrival, while in $\mathscr{Q}_3$ it starts
with delay, resulting in \eqref{3.1}.
\end{proof}

\begin{prop}\label{prop2}
\begin{equation}\label{3.2}
Q_3(t,\omega)-Q_1(t,\omega)\leq1.
\end{equation}
\end{prop}

\begin{proof}
 For the purpose of the
proof we will follow up the sample paths of the both processes of
the queueing systems $\mathscr{Q}_1$ and $\mathscr{Q}_3$.

Apparently, that $Q_1(t,\omega)=Q_3(t,\omega)=0$ for all
$t\in[0,t_1)$ (recall that according to convention $A(0)=D(0)=0$).
 For the system $\mathscr{Q}_1$, let
 $$l=\inf\left\{i: x_i>t_1\right\}.$$ Then, for any $t$ from the
 interval $[t_1, x_l)$ we have $Q_1(t,\omega)=Q_3(t,\omega)=A(t,\omega)$,
 while in the point $x_l$ itself we have
 $Q_3(x_l, \omega)-Q_1(x_l,\omega)=1$, see Figure 1 (a).
\begin{figure}
\includegraphics[width=15cm,height=20cm]{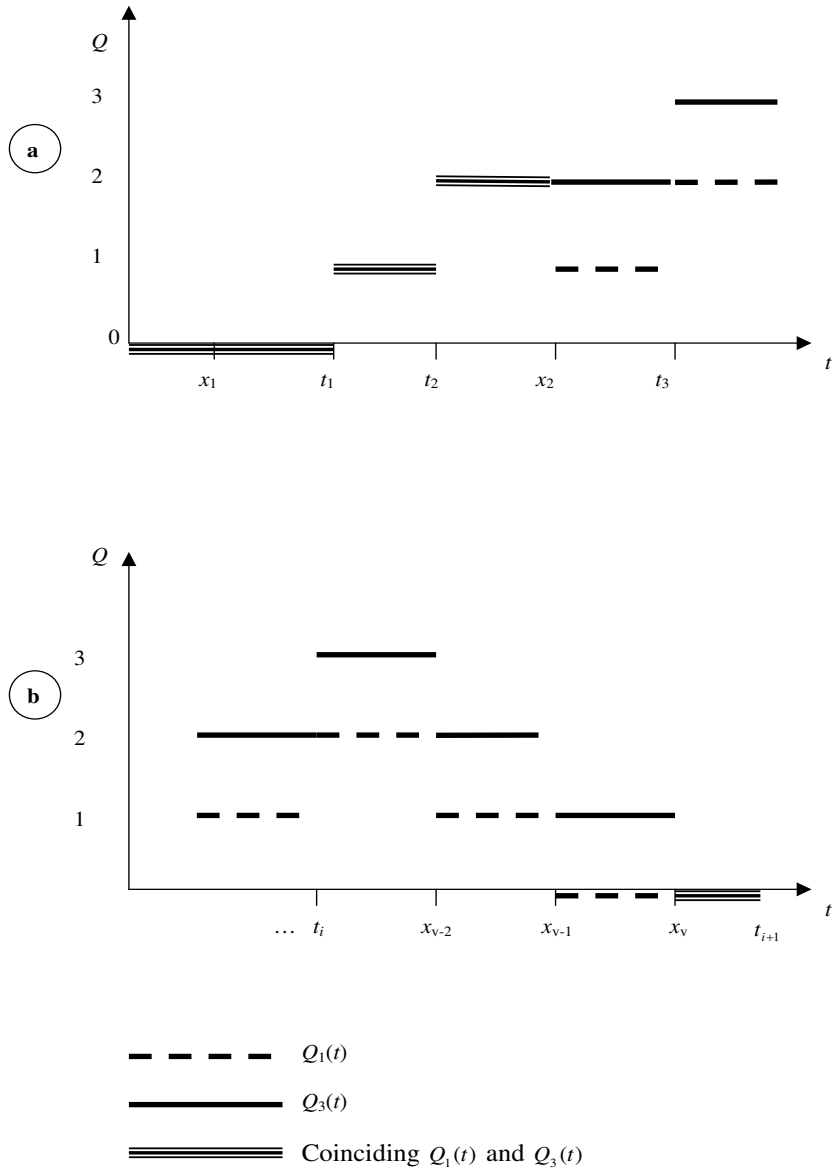}
\caption{Typical sample paths of the queue-length processes
$Q_1(t)$ and $Q_3(t)$:\newline (a) Sample paths of the
queue-length processes after the origin;\newline (b) Sample paths
of the queue-length processes before approaching zero at point
$x_v$.}
\end{figure}

 Let $v$ be a new number greater
 than $l$ satisfying the property
 $$v=\inf\left\{k>l: Q_1(x_k-)=0\right\}.$$
(If such a number does not exist, then $v$ equates to infinity. In
this case, obviously, $Q_3(t, \omega)-Q_1(t,\omega)=1$, $t\geq
x_l$.) Then for all $t$ of the interval $[x_l, x_v)$ we have
 $Q_3(t,\omega)-Q_1(t,\omega)=1$, while in the point $x_v$ itself we arrive at
 $Q_1(x_v,\omega)=Q_3(x_v,\omega)=0$, see Figure 1 (b).

 Thus, we arrived at zeroth queue-lengths again.
 The further paths of the both processes after point
$t=x_v$ behave similarly to those after the point $t=0$, i.e. the
difference $Q_3(t,\omega)-Q_1(t,\omega)$ can take only one of the
two values 0 or 1.
\end{proof}

From Propositions \ref{prop1} and \ref{prop2} we arrive at the
following conclusion: \textit{if the queueing system
$\mathscr{Q}_1$ is stable, then the both queueing systems
$\mathscr{Q}_2$ and $\mathscr{Q}_3$ are stable as well}. Following
Proposition \ref{prop1} this statement of stability can be
extended for more complicated constructions including two, three
and more different arrival processes and can be then applied to
queueing networks.

\section{Stability of JS-queues}

In this section we study stability of JS-queues. We start from the
simplest case of symmetric $\Sigma_m$ queues. Denote
$D_a(t)=\sum_{j=1}^m D_a^{(j)}(t)$. In Theorem \ref{thm1} below,
the assumption that the processes $A(t)$, $A^\prime(t)$ and
$D^{(j)}(t)$, $j=1,2,\ldots,m$, all are mutually independent point
processes is relaxed. For $\Sigma_m$ queues where
$p_j=\frac{1}{m}$, $j=1,2,\ldots,m$, instead of \eqref{r1.1},
\eqref{r1.2} and \eqref{r1.3} we assume
\begin{equation}\label{4.0}
\mathbf{P}\left\{\lim_{a\to~-\infty}\frac{A_a(t)+A_a^{\prime}(t)-D_a(t)}{t-a}
={\lambda+\lambda^\prime-\mu m}\right\}=1.
\end{equation}
%

\begin{thm}
\label{thm1} In addition to \eqref{4.0} assume that
\begin{equation}
\label{4.1}
\lim_{a\to~-\infty}\mathbf{P}\{A_a(t)+A_a^\prime(t)-D_a(t)\in\mathscr{S}\}=0
\end{equation}
for any bounded set $\mathscr{S}$. Then, the system $\Sigma_m$ is
stable if and only if the condition
$\frac{\lambda}{m}+\frac{\lambda^\prime}{m}<\mu$ is fulfilled.
\end{thm}

\begin{proof}
Since the families $\{A_a^{(j)}(t)\}_{j\leq m}$
 and $\{D_a^{(j)}(t)\}_{j\leq m}$
consist of identically distributed processes, then the family
$\{A_a^{\prime (j)}(t)\}_{j\leq m}$ also consists of identically
distributed processes, and according to (\ref{r1.4}) the family of
the processes
 $\{Q_a^{(j)}
(t)\}_{j\le m}$ consists of identically distributed processes too.
Observe that from (\ref{4.1}), because the system $\Sigma_m$ is
symmetric, we also have
\begin{equation}\label{4.1'}
\lim_{a\to~-\infty}\mathbf{P}\{A_a^{(j)}%
(t)+A_a^{\prime(j)}(t)-D_a^{(j)}(t)\in\mathscr{S}\}=0
\end{equation}
for all $j=1,2,\ldots,m$.

Notice, that if
$\frac{\lambda}{m}+\frac{\lambda^{\prime}}{m}<\mu$, then
\begin{equation}\label{4.2}
\mathbf{P}\Big\{\sup_{a\leq t<\infty}\big[A_a^{(j)}(t)+A_a^{\prime(j)}%
(t)-D_a^{(j)}(t)]<\infty\Big\}=1.
\end{equation}
Indeed, according to \eqref{4.0}, $\mathbf{P}$-a.s.
\begin{equation}\label{4.3}
\lim_{a\to~-\infty}\frac{A_a^{(j)}(t)+A_a^{\prime(j)}(t)-D_a^{(j)}(t)}{t-a}%
=\frac{\lambda}{m}+\frac{\lambda^{\prime}}{m}-\mu<0.
\end{equation}
Hence, $\mathbf{P}$-a.s.
\[
\lim_{a\to~-\infty}[A_a^{(j)}(t)+A_a^{\prime(j)}(t)-D_a^{(j)}(t)]=-\infty,
\]
and (\ref{4.2}) follows from the fact that $A_a^{(j)}(t)$,
$A_a^{\prime(j)}(t)$ and $D_a^{(j)}(t)$ all are c\'adl\'ag
processes.

Next, taking into account \eqref{4.3} and
$A_a^{(j)}(a)=A_a^{\prime(j)}(a)=D_a^{(j)}(a)=0$, for the
queue-length process $Q^{(j)}(t)$ one can write the following
representation:
\begin{equation}\label{4.4}
Q_a^{(j)}(t)=\big[A_a^{(j)}(t)+A_a^{\prime(j)}(t)-D_a^{(j)}(t)\big]-\inf_{a\leq
s\leq t}\big[A_a^{(j)}(t)+A_a^{\prime(j)}(t)-D_a^{(j)}(t)\big].
\end{equation}
Representation \eqref{4.4} is well-known (e.g. Borovkov
\cite{Borovkov (1976)}) and is a consequence from the Skorokhod
reflection principle (e.g. Kogan and Liptser \cite{Kogan and
Liptser (1993)} for typical application to queue-length
processes).

Next, from \eqref{4.4}
we have:
\begin{eqnarray}\label{4.5}
Q_a^{(j)}(t)&{\buildrel d\over =}&\sup_{a\leq s\leq
t}\big\{\big[A_a^{(j)}(t)+A_a^{\prime(j)}
(t)-D_a^{(j)}(t)\big]\nonumber\\
&&-\big[A_a^{(j)}(s)+A_a^{\prime(j)}(s)-D_a^{(j)}(s)\big]\}.
\end{eqnarray}
From representation \eqref{4.5}, relation \eqref{4.2} and the fact
that all of the processes $A_a^{(j)}(t)$, $A_a^{\prime(j)}(t)$ and
$D_a^{(j)}(t)$ are c\'adl\'ag processes, it follows that there
exists a bounded set $\mathscr{S}_0$ such that
\[
\limsup_{a\to~-\infty}\mathbf{P}\{Q^{(j)}(t)\in\mathscr{S}_0\}>0,
\]
and the sufficient condition
is therefore proved. The necessary condition follows from the fact
that (\ref{4.1'}) together with (\ref{4.2}) imply the condition
$\frac{\lambda}{m}+\frac{\lambda^{\prime}}{m}<\mu$.
\end{proof}

\begin{rem}\label{rem0}
If $\frac{\lambda}{m}+\frac{\lambda^{\prime}}{m}=\mu$, and the
processes $A(t)$, $A^\prime(t)$ and $D(t)$ all are non-trivial
renewal processes, then we easily arrive at condition \eqref{4.1}.
However, there are examples where
$\frac{\lambda}{m}+\frac{\lambda^{\prime}}{m}=\mu$, but condition
\eqref{4.1} is not fulfilled. Indeed, let
$\tau_1+\tau_1^\prime-\chi_1$ be a uniformly distributed random
variable in $[-b, b]$, ($b>0$), and let
$\tau_{i+1}+\tau_{i+1}^\prime-\chi_{i+1}=-(\tau_{i}+\tau_{i}^\prime-\chi_{i})$,
$i\geq1$. In this case \eqref{4.1} is not valid. Therefore
condition \eqref{4.1} is meaningful.
\end{rem}
\begin{rem}The conditions of Theorem \ref{thm1} are weaker than
conditions of other known stability theorems (e.g. Borovkov
\cite{Borovkov (1976)}) requiring the strict stationarity and
ergodicity of appropriate sequence of random variables. For
example, assume that the sequence
$\{\tau_n+\tau_n^\prime-\chi_n\}_{n\geq1}$ consists of identically
distributed random variables, $\tau_1+\tau_1^\prime-\chi_1$ is a
uniformly distributed random variable in $[-b, b]$, ($b>0$), and
$\tau_2+\tau_2^\prime-\chi_2=\tau_1+\tau_1^\prime-\chi_1$, and
$\tau_{i+1}+\tau_{i+1}^\prime-\chi_{i+1}=-(\tau_{i}+\tau_{i}^\prime-\chi_{i})$,
$i\geq2$. Then the above sequence
$\{\tau_n+\tau_n^\prime-\chi_n\}_{n\geq1}$ is not strictly
stationary. However, the queueing system associated with this
sequence is stable.
\end{rem}

For following Theorem \ref{thm2}, our assumptions regarding the
point processes $A_a(t)$, $D_a^{(j)}(t)$ is as follows. For
dependent processes $A_a(t)$, $D_a^{(j)}(t)$ we suppose that the
normalized processes $\frac{A_a(t)}{t-a}$ and
$\frac{D_a^{(j)}(t)}{t-a}$ converge, as $a\to~-\infty$, to the
corresponding limits $\lambda$ and $\mu$ only in distribution,
while
\begin{equation}\label{4.8}
\mathbf{P}\left\{\lim_{a\to~-\infty}\frac{A_a(t)-D_a^{(j)}(t)}{t-a}=\lambda-\mu\right\}=1.
\end{equation}
Then, the process $A_a^\prime(t)$ is assumed to satisfy the
condition
\begin{equation}\label{4.9}
\mathbf{P}\left\{\lim_{a\to~-\infty}\frac{A_a^\prime(t)}{t-a}=\lambda^\prime\right\}=1.
\end{equation}

Next Theorem \ref{thm2} is related to stability of $\wp_{m}$
queues. Denote $\lambda_j=\lambda p_j$, $j^*=\arg\max_{1\leq j\leq
m}\lambda_j$ and $\Delta=\sum_{j=1}^m(\lambda_{j^*}-\lambda_j)$.

\begin{thm}
\label{thm2} In addition to \eqref{4.8} and \eqref{4.9} assume
that both
\begin{equation}\label{4.10}
\lim_{a\to~-\infty}\mathbf{P}\{A_a(t)+A_a^\prime(t)-D_a(t)\in\mathscr{S}\}=0,
\end{equation}
and
\begin{equation}\label{4.11}
\lim_{a\to~-\infty}\mathbf{P}\{A_a^{(j^*)}(t)-D_a^{(j^*)}(t)\in\mathscr{S}\}=0
\end{equation}
for any bounded set $\mathscr{S}$. Then the system is stable if
and only if one of the following two conditions is fulfilled:
$$
\begin{cases}
\lambda^{(j^*)}<\mu, &\mbox{if} \ \Delta\geq\lambda^\prime,\\
\lambda+\lambda^\prime<m\mu, &\mbox{otherwise}.
\end{cases}
$$
\end{thm}

\begin{proof}
The theorem can be proved by a slight modification of the earlier
proof.

Suppose first that $\Delta\ge\lambda^{\prime}$, and denote $q_j$
the fraction of customer of the opportunistic traffic being
assigned to the $j$th queue. From the limiting relations
\begin{eqnarray*}
\lim_{a\to~-\infty}\frac{A_a^{(j)}(t)}{t-a}&{\buildrel d\over
=}&\lambda_j,\\
\lim_{a\to~-\infty}\frac{A_a^{\prime(j)}(t)}{t-a}&{\buildrel
{a.s.}\over =}&\lambda_j^\prime,
\end{eqnarray*}
according to well-known Skorokhod's theorem \cite{Skorokhod
(1956)}, p. 281, one can conclude that there exists a probability
space $(\Omega, \mathscr{F}, \mathbf{P})$, with given there a
family of processes
$\left\{\frac{A_a^{(j)}(t,\omega)+A_a^{\prime(j)}(t,\omega)}{t-a},
t>a\right\}$ such that for $\mathbf{P}-$ almost all
$\omega\in\Omega$,
\begin{eqnarray*}
\lim_{a\to~-\infty}\frac{A_a^{(j)}(t,\omega)+A_a^{\prime(j)}(t,\omega)}{t-a}&=&\lambda_j+\lambda_j^\prime.
\end{eqnarray*}
Therefore, from the balance equations
\[
\lambda_j+\lambda_j^\prime=\lambda_j+\lambda^\prime q_j=\varrho,
\]
one can conclude that $q_{j^*}$ must be equal to 0. Therefore
\[
\lambda_{j^*}\geq \lambda_j+\lambda_j^\prime
\]
for all $j=1,2,\ldots,m$, where $\lambda_j^\prime=\lambda^\prime
q_j$. Then, in this probability space for $\mathbf{P}-$ almost all
$\omega\in\Omega$,
\[
\lim_{a\to~-\infty}\frac{A_a^{(j)}(t,\omega)+A_a^{\prime(j)}(t,\omega)}{t-a}\leq
\lim_{a\to~-\infty}\frac{A_a^{(j^*)}(t,\omega)}{t-a},
\]
and, from \eqref{r1.5} and coupling arguments of sample paths
comparison
\[
\lim_{a\to~-\infty}\frac{Q_a^{(j)}(t,\omega)}{t-a}\leq
\lim_{a\to~-\infty}\frac{Q_a^{(j^*)}(t,\omega)}{t-a}.
\]

Therefore, in the original probability space we have
\begin{equation*}\label{4.12}
\lim_{a\to~-\infty}\frac{A_a^{(j)}(t)+A_a^{\prime(j)}(t)}{t-a}~{\buildrel
d\over\leq}~ \lim_{a\to~-\infty}\frac{A_a^{(j^*)}(t)}{t-a},
\end{equation*}
and consequently,
\begin{equation*}
\lim_{a\to~-\infty}\frac{Q_a^{(j)}(t)}{t-a}~{\buildrel
d\over\leq}~ \lim_{a\to~-\infty}\frac{Q_a^{(j^*)}(t)}{t-a}
\end{equation*}
for all $j=1,2,\ldots,m$. Hence the problem reduces to conditions
of stability of a single queueing system with autonomous service
mechanism with the given arrival process $A^{(j^*)}(t)$ and
departure process $D^{(j^*)}(t)$, and under assumptions
\eqref{4.8} and \eqref{4.11} the necessary and sufficient
condition of stability is given by $\lambda_{j^*}<\mu$.

Let us now consider the opposite case $\Delta<\lambda^\prime$ and
assumption $\lambda+\lambda^\prime<\mu m$. Then, there exist
probabilities $q_j$, $j=1,2,\ldots,m$, all strictly positive, and
$\sum_{j=1}^m q_j=1$. Under these probabilities, the opportunistic
traffic is thinned into $m$ processes such that almost surely
\[
\lim_{a\to~-\infty}\frac{A_a^{\prime(j)}(t)}{t-a}=\lambda^\prime
q_j=\lambda_j^\prime.
\]

Indeed, since $\lambda+\lambda^\prime<\mu m$, then there exists
the value $\varrho=\frac{\lambda+\lambda^\prime}{m}<\mu$ such that
for all $j$
\[
\varrho>\lambda_j,
\]
and therefore,
\begin{equation}\label{4.13}
q_j=\frac{\varrho-\lambda_j}{\lambda^\prime}>0.
\end{equation}
Since $\lambda_j+q_j\lambda^\prime=\varrho$, the family
$\left\{A_a^{(j)}(t)+A_a^{\prime(j)}(t)\right\}_{j\leq m}$
consists of processes having the same rate, i.e.
\begin{equation}
\label{4.14}
\lim_{a\to~-\infty}\frac{A_a^{(j)}(t)+A_a^{\prime(j)}(t)}{t-a}~{\buildrel
d\over =}~\varrho.
\end{equation}
The possible values $\{q_{j}\}_{j\leq m}$ are unique, since
otherwise, if there are different arrival rates
$\lambda_{j}+\lambda^{\prime}q_{j}^\prime$, then one of queues
must be stochastically longer than other. Let $j^{*}$ be the order
number of the longer queue. Then $q_{j^{*}}$ should be equal to 0,
and we have the contradiction with \eqref{4.13}.

According to (\ref{4.14}) for each of $m$ queue-length processes
the arrival rate is the same. With the same arrival and departure
intensities one can repeat the proof of Theorem \ref{thm1} for
each queue-length process. Therefore, $\varrho<\mu$ is the
condition for stability, and under condition \eqref{4.10} the
system is stable if and only if $\varrho<\mu $. The theorem is
proved.
\end{proof}

\begin{rem}
\label{rem2} In the case of the model $\Gamma_m$ the following
additional comments are necessary. If the value $j^*$ is not
unique, then condition \eqref{4.11} should be assumed for all
$j^*$. In addition, instead of condition \eqref{4.8} we should
require
\[
\mathbf{P}\left\{\lim_{a\to~-\infty}\frac{A_a^{(j)}(t)-D_a^{(j)}(t)}{t-a}=\lambda_j-\mu\right\}=1
\]
for all $j=1,2,\ldots,m$, and in addition we should require the
convergence of the normalized processes $\frac{A_a^{(j)}(t)}{t-a}$
and $\frac{D_a^{(j)}(t)}{t-a}$, as $a\to~-\infty$, to the
corresponding numbers $\lambda_j$ and $\mu$ in distribution.
\end{rem}

\section{Load-balanced networks and their stability}
In the previous section, the necessary and sufficient conditions
for stability of systems $\wp_{m}$ have been established. In this
section we extend above Theorem \ref{thm2} for load-balanced
networks associated with $\wp_{m}$ queueing systems. The main
result of this section is Theorem \ref{thm5} establishing the
necessary and sufficient condition for the stability of
load-balanced networks. Theorems \ref{thm3} and \ref{thm4} are the
preliminary results establishing only sufficient conditions for
the stability.

The load-balanced network considered below is the following
extension of the $\wp_{m}$ queueing system.

Assume that an arriving customer of the dedicated traffic occupies
the server $j$ with probability $p_j$ \ ($j=1,2,\ldots,m$), and
$\sum_{j=1}^{m}p_{j}=1$. After his service completion in the $j$th
queue, a customer leaves the system with probability
$1-p_{j}^{*}$, remains at the same $j$th queue with probability
$p_{j,j}$, goes to the different queue $i\neq j$ with probability
$p_{j,i}$, and choose the shortest queue with probability
$p_{j,sh}$, breaking ties at random. It is assumed that
$p_{j}^{*}<1$ at least for one of the indexes $j$. This model is
called \textit{load-balanced network}. The stability conditions of
the Markovian variant of this network containing two stations only
has been established by Kurkova \cite{Kurkova (2001)}.

Some different variants of this network have also been studied in
Martin and Suhov \cite{Martin and Suhov (1999)}, Vvedenskaya
Dobrushin and Karpelevich \cite{Vvedenskaya Dobrushin and
Karpelevich (1996)} and other papers.

Denote $\lambda_{j}=\lambda p_{j}$, $
\Lambda_{j}=\lambda_{j}+\mu\sum_{i=1}^{m}p_{i,j}$,
$j^{*}=\arg\max_{1\le j\le m}\Lambda_{j}$.

For the sake of simplification of the proofs, we follow the
assumptions that are described in Section 2, with understanding
that they can be relaxed by the way of the previous section. We
also assume that the processes all are started at $a$ and use the
variants of the assumptions where $a\to~-\infty$ rather than
$t\to\infty$.

The sufficient condition given in Theorem \ref{thm3} is a
straightforward extension  of earlier Theorem \ref{thm2} written
now in a simpler form.

\begin{thm}\label{thm3}
Assume that $\lambda_{j^{*}}\geq\Lambda_{j}$ for all $j\ne j^{*}$.
Denote
$$\Delta _{1}=\sum_{j=1}^{m}(\Lambda_{j^{*}}-\lambda_{j}),$$
and
$$\Delta_{2}=\sum_{j\neq
j^{*}}(\lambda_{j^{*}}-\Lambda_{j}).$$
Then the load-balanced
network is stable if one of the following two conditions is
fulfilled:
\[%
\begin{cases}%
\Lambda_{j^{*}}<\mu, &\mbox{if} \ \Delta_{2}
\geq\lambda^{\prime}+\mu
\sum_{j=1}^{m}p_{j,sh},\\
\lambda+\lambda^{\prime}+\mu\Big(\sum_{i=1}^{m}\sum_{j=1}^{m}p_{i,j}%
+\sum_{j=1}^{m}p_{j,sh}\Big)<m\mu, &\mbox{if} \
\Delta_{1}<\lambda^{\prime}.
\end{cases}
\]
\end{thm}

\begin{proof} The proof of the theorem starts from the case
$p_{j,sh}=0$ for all $j=1,2,\ldots,m$ and then discusses the case
$p_{j,sh}\geq 0$ for all $j=1,2,\ldots,m$.

Let us start from the case where $p_{j,sh}=0$ for all
$j=1,2,\ldots,m$.
Then, the system of equations for the
queue-length processes started at $t=a$ can be written
\begin{equation}\label{5.1}
Q_{a}^{(j)}(t)=A_{a}^{(j)}(t)+A_{a}^{\prime(j)}(t)
+E_{a}^{(j)}(t)-\int_{a}^{t}\mathbf{1}%
_{\{Q_{a}^{(j)}(s-)>0\}}\mbox{d}D_{a}^{(j)}(s),
\end{equation}
where the new process $E_{a}^{(j)}(t)$ presenting in \eqref{5.1}
is a process, generated by \textit{internal dedicated arrivals} to
the $j$th queue. By internal dedicated arrivals to the $j$th queue
we mean internal arrivals of the customers, who after their
service completion in one or other queue are assigned to the $j$th
queue. Recall that there is probability $p_{i,j}$ to be assigned
from the queue $i$ to the queue $j$. Relationship \eqref{5.1} is
of the same type as that \eqref{r1.5}, and therefore all of the
arguments of the earlier proof of Theorem \ref{thm2} can be
repeated. Specifically, in the case $\Delta_2\geq\lambda^{\prime}$
the system is stable if $\Lambda_{j^{*}}<\mu $.

In
turn, the process $E_{a}^{(j)}(t)$ can be represented as $\sum_{i=1}^{m}%
E_{a}^{(i,j)}(t)$, where the point process $E_{a}^{(i,j)}(t)$ is
generated by the customers who are assigned to the $j$th queue
after their service completion in the $i$th queue (in the case
$i=j$ it is assumed that the customers decide to stay at the same
queue). Notice, that
$$
\lim_{a\to~-\infty}\frac{E_a^{(i,j)}(t)}{t-a}
$$
does exist with probability 1 (because the process
$E_a^{(i,j)}(t)$ is generated by the procedure of thinning of the
departure process $D_a^{(i)}(t)$), and
\begin{equation}\label{5.2}
\mathbf{P}\Big\{\lim_{a\to~-\infty}\frac{E_a^{(i,j)}(t)}{t-a}\leq\mu
p_{i,j}\Big \}=1,
\end{equation}
where the equality holds only in the case where the fraction of
the $i$th queue idle period vanishes as $a\to~-\infty$. Therefore
under the condition $\Delta_1<\lambda^{\prime}$ the system is
stable if $\lambda+\lambda^\prime+\mu\sum_{i=1}^m\sum_{j=1}^m
p_{i,j}<m\mu$.

Assume now that $p_{j,sh}\geq0$ for all $j=1,2,\ldots,m$. Then
instead of \eqref{5.1} we have the equation
\begin{equation}\label{5.3}
\begin{aligned}
Q_{a}^{(j)}(t)&=A_{a}^{(j)}(t)+A_{a}^{\prime(j)}(t)
+E_{a}^{(j)}(t)+E_{a}^{\prime(j)}(t)\\
&\ \ \ -\int_{a}^{t}\mathbf{1}
_{\{Q_{a}^{(j)}(s-)>0\}}\mbox{d}D_{a}^{(j)}(s),
\end{aligned}
\end{equation}
where $E^{\prime(j)}(t)$ are the point processes associated with
\textit{internal opportunistic traffic} to the $j$th queue. By
internal opportunistic traffic we mean the internal traffic of
customers presenting in the queue, who after their service
completion decide to join the shortest queue. In the case where
the shortest queue is the $j$th queue, we say about opportunistic
traffic to the $j$th queue. Again,
$$
\lim_{a\to~-\infty}\frac{E_a^{\prime(j)}(t)}{t-a}
$$
does exist with probability 1 (because the process
$E_a^{\prime(j)}(t)$ is generated by the procedure of thinning of
the departure process $D_a^{(j)}(t)$), and similarly to
\eqref{5.2} we have:
\begin{equation}\label{5.4}
\mathbf{P}\left\{\lim_{a\to~-\infty}\frac{\sum_{j=1}^{m}E_a^{\prime(j)}(t)}%
{t-a}\leq\mu\sum_{j=1}^{m}p_{j,sh}\right\}=1.
\end{equation}
Therefore, the entire opportunistic traffic to the $j$th queue,
being a sum of the processes $A_a^{\prime(j)}(t)$ and
$E_a^{\prime(j)}(t)$, satisfies
\begin{equation}\label{5.5}
\mathbf{P}\left\{\lim_{a\to~-\infty}\frac{\sum_{j=1}^{m}[A_a^{\prime%
(j)}(t)+E_a^{\prime(j)}(t)]}{t-a}\leq\lambda^{\prime}+\mu\sum_{j=1}^{m}%
p_{j,sh}\right\}=1.
\end{equation}
Hence the proof of this theorem remains similar to the proof of
Theorem \ref{thm2}. In the case $\Delta_{2}
\geq\lambda^{\prime}+\mu \sum_{j=1}^{m}p_{j,sh}$ the system is
stable if $\Lambda_{j^{*}}<\mu $. In the other case
$\Delta_{1}<\lambda^{\prime}$ the system is stable if
$\lambda+\lambda^{\prime}+\mu\Big(\sum_{i=1}^{m}\sum_{j=1}^{m}p_{i,j}%
+\sum_{j=1}^{m}p_{j,sh}\Big)<m\mu$. The conditions of the theorem
are sufficient and not necessary, because the left-hand sides of
\eqref{5.2}, \eqref{5.4} and \eqref{5.5} contain the probability
of inequalities, and the exact parameters of internal dedicated
traffic as well as internal opportunistic traffic are unknown.
\end{proof}

In order to formulate and prove a necessary and sufficient
condition of stability for the above load-balanced network, we
first need to improve the sufficient condition given by Theorem
\ref{thm3}. For this purpose, rewrite \eqref{5.2} as
\[
\mathbf{P}\Big\{\lim_{a\to~-\infty}\frac{E_a^{(i,j)}(t)}{t-a}=\varrho_{i}\mu
p_{i,j}\Big\}=1,
\]
where the value $\varrho_{i}$ satisfies the inequality
$0<\varrho_{i}\leq 1$. The value $\varrho_{i}$ is the fraction of
time when the server of the $i$th queue is busy. Then the case of
$\varrho_{i}=1$ means that the server of the $i$th queue is busy
almost always.

\smallskip
Let us consider the system of inequalities
\begin{equation}\label{5.6}
\frac{\lambda_{j}+\mu\sum_{i=1}^{m}\varrho_{i}p_{i,j}}{\varrho_{j}^{*}}\leq\mu,
\ j=1,2,\ldots,m.
\end{equation}
The meaning of inequality \eqref{5.6} is the following. The
left-hand side contains the total sum of rates of dedicated
traffic to the $j$th queue divided to the traffic parameter
$\varrho_{j}^{*}$ of the $j$th queue. The total sum of rates of
dedicated traffic of the $j$th queue consists of exogenous and
internal arrivals to that $j$th queue, excluding the rates for
joining the shortest queue customers. Since the rates associated
with opportunistic traffic are excluded, there is the inequality
'$\le$' between the left and right sides. Thus, if the $j$th queue
is never shortest, then the sum of the rates of the left-hand side
divided to $\varrho_{j}^{*}$ becomes equal to $\mu$ of the
right-hand side. When the traffic parameter $\varrho_{j}^{*}$ is
greater than 1, the $j$th queue increases to infinity with
probability 1.
Therefore, in the sequel we only consider the case when $\varrho_{j}^{*}%
\leq 1$ for all $j=1,2,\ldots,m$. In this case
$\varrho_{j}^{*}=\varrho_{j}$, and we therefore have
\begin{equation}\label{5.7}
\lambda_{j}+\mu\sum_{i=1}^{m}\varrho_{i}p_{i,j}\le\varrho_{j}\mu,
\ j=1,2,\ldots,m.
\end{equation}

Let us now write a so-called \textit{balance equation}, taking
into account also joining the shortest queue customers. We have
\begin{equation}\label{5.8}
\lambda+\lambda^{\prime}+\mu\sum_{j=1}^{m}\sum_{i=1}^{m}\varrho_{i}p_{i,j}
=\mu\sum_{j=1}^{m}\varrho_{j}(1-p_{j,sh}).
\end{equation}
Now, we are ready to prove the \textit{improved} sufficient
condition for the stability. This version is also based on a
straightforward extension of Theorem \ref{thm2}.

\begin{thm}\label{thm4}
The load-balanced network is stable if there exists
\[
\varrho^{*}=\max_{1\le j\le m}\varrho_{j},
\]
satisfying the condition $\varrho^*<1$, where the values
$\varrho_j$, $j=1,2,\ldots,m$, are defined by \eqref{5.7} and
\eqref{5.8}.
\end{thm}

\begin{proof}
Let $\Lambda_{j}=\lambda_{j}+\mu\sum_{i=1}^{m}%
\varrho_{i}p_{i,j}$, let $j^{*}=\arg\max_{1\le j\le
m}\Lambda_{j}$, and let
$\Delta=\sum_{j=1}^{m}(\Lambda_{j^{*}}-\Lambda_{j})$. In the case
$\Delta
>\lambda^{\prime}+\mu\sum_{j=1}^{m}\varrho_{j}p_{j,sh}$ the $j^{*}$-th
queue is never shortest, and therefore, following the proof of
Theorem \ref{thm2}, the system is stable if $\Lambda_{j^{*}}<\mu$.
Therefore from \eqref{5.6} we have $\varrho_{j^{*}}<1$, and
because the $j^{*}$-th queue is the longest queue, we have
$\varrho^*=\varrho_{j^{*}}\ge\varrho_{j}$, $j=1,2,\ldots,m$.
Therefore $\varrho_{j}<1$ for all $j=1,2,\ldots,m$ is a sufficient
condition of stability for this case.

Let us now consider the opposite case, where $\Delta\leq\lambda^{\prime}%
+\mu\sum_{j=1}^{m}\varrho_{j}p_{j,sh}$. As in the proof of Theorem
\ref{thm2}, in this case the arrival rate to all $m$ queues is the
same, and therefore $\varrho_{1}=\varrho_{2}=\ldots=\varrho_{m}$.
Thus, the only two cases are there as $\varrho_{j}<1$ or
$\varrho_{j}=1$ for all $j=1,2,\ldots ,m$. In the case
$\varrho_{j}<1$, the stability result is analogous to that of
Theorem \ref{thm2}, since in this case from \eqref{5.8} we obtain
\[
\lambda+\lambda^{\prime}+\mu\sum_{j=1}^{m}\sum_{i=1}^{m}p_{i,j}+\mu\sum
_{j=1}^{m}p_{j,sh} <m\mu.
\]
The theorem is proved.
\end{proof}

Now in order to formulate and prove a necessary and sufficient
condition for stability, let us consider the following linear
programming problem in $\mathbb{R}^{m+1}$:
\begin{equation}\label{5.10}
\mbox{Minimize} \ x_{m+1}
\end{equation}
subject to the restrictions:
\begin{equation}\label{5.11}
\lambda_{j}+\mu\sum_{i=1}^{m}x_{i}p_{i,j}\le x_{j}\mu, \
j=1,2,\ldots ,m,
\end{equation}%
\begin{equation}\label{5.12}
\lambda+\lambda^{\prime}+\mu\sum_{j=1}^{m}\sum_{i=1}^{m}x_{i}p_{i,j}
=\mu \sum_{j=1}^{m} x_{j}(1-p_{j,sh}),
\end{equation}%
\begin{equation}\label{5.13}
x_{j}\le x_{m+1}, \ j=1,2,\ldots,m.
\end{equation}

Observe, that the restrictions \eqref{5.11} and \eqref{5.12}
correspond to \eqref{5.7} and \eqref{5.8}, where the values
$\varrho_{j}$ are replaced with unknown $x_{j}$. The functional of
\eqref{5.10} and inequalities \eqref{5.13} are associated with the
condition of Theorem \ref{thm4}: $\max_{1\leq j\leq
m}\varrho_{j}<1$. $x_{m+1}$ is an additional variable; thus the
linear programming \eqref{5.10}-\eqref{5.13} is a mini-max
problem. That is, if the minimum of $x_{m+1}$ is achieved in some
point $x_{m+1}^{*}<1$, then all components of the vector
($x_{1}^{*}, x_{2}^{*},\ldots,x_{m+1}^{*}$) associated with this
solution are less than 1, and there exists a solution of system
\eqref{5.7} and \eqref{5.8} with $\varrho_{j}<1$ for all
$j=1,2,\ldots,m$. Therefore in the following the vector associated
with solution of the problem \eqref{5.10}-\eqref{5.13} is denoted
($\varrho_{1}, \varrho_{2},\ldots ,\varrho_{m}$). Otherwise if
$x_{m+1}^{*}\geq 1$, then we set $\varrho_{j}=1$,
$j=1,2,\ldots,m$.

\smallskip
Next, denote $A_a^{(j)}(t)+E_a^{(j)}(t)$ the dedicated arrival
process. Its relation to the initial processes $A_a(t)$ and
$D_a^{(j)}(t)$ is the following. Each arrival of the initial
process $A_a(t)$ is forwarded to the queue $j$ with probability
$p_{j}$, and each customer served in the queue $i$ returns to the
$j$th queue with
probability $p_{i,j}$. Then, the process $A_a^{(j)}%
(t)+E_a^{(j)}(t)$ is a sum of all arrivals of external and
internal dedicated traffic, and
\[
\mathbf{P}\Big\{\lim_{a\to~-\infty}\frac{A_a^{(j)}(t)+E_a^{(j)}(t)}{t-a}=\Lambda
_{i}=\lambda_{j} +\mu\sum_{i=1}^{m}\varrho_{i}p_{i,j}\Big\}=1,
\]
where $\varrho_{j}$, $j=1,2,\ldots,m$, are a solution of the
linear programming given by \eqref{5.10}-\eqref{5.13}. Now let
$j^{*}=\arg\max_{1\le j\le m}\Lambda_{j}$. We have the following
theorem.

\begin{thm}\label{thm5}
Assume that the both
\begin{equation}\label{5.15}
\lim_{a\to~-\infty}\mathbf{P}\{A_a^{(j^{*})}(t)+E_a^{(j^{*})}(t)-D_a^{(j^{*})}%
(t)\in\mathscr{S}\}=0,
\end{equation}
and
\begin{equation}\label{5.16}
\lim_{a\to~-\infty}\mathbf{P}\{A_a(t)+A_a^{\prime}(t)+E_a(t)+E_a^{\prime}(t)-D_a(t)\in
\mathscr{S}\}=0,
\end{equation}
for any bounded set $\mathscr{S}$, where
$E_a(t)=\sum_{j=1}^{m}E_a^{(j)}(t)$ is the point process
associated with all internal arrivals of dedicated traffic,
$E_a^{\prime}(t)$ is the point process associated with all
internal arrivals of opportunistic traffic,
$D_a(t)=\sum_{j=1}^{m}D^{(j)}(t)$.  Then the load-balanced network
is stable if and only if $\max_{1\leq j\le m}\varrho_{j}<1$.
\end{thm}

\begin{rem}\label{rem5-1}
Conditions \eqref{5.15} and \eqref{5.16} are verifiable
conditions. As soon as the linear programming problem is solved
and we know the vector of solution $(\varrho_1,
\varrho_2,\ldots,\varrho_m)$, the unknown processes $E_a(t)$ and
$E_a^\prime(t)$ as well as $E_a^{(j)}(t)$ and
$E_a^\prime{}^{(j)}(t)$ can be easily modelled via derivative
processes $D_a^{(j)}(t)$ ($j=1,2,\ldots,m$). Note also, that above
conditions \eqref{5.15} and \eqref{5.16} are automatically
fulfilled if $\max_{1\leq j\leq n}\varrho_j<1$.
\end{rem}
\begin{rem}\label{rem5-2}
In the case of the network associated with the model $\Gamma_m$
the following additional comment is necessary. If the value
$j^{*}$ is not unique, then condition \eqref{5.15} should be
assumed for all of these values $j^{*}$.
\end{rem}

\section{Concluding remarks}
In this paper we established stability of different type
joint-the-shortest-queue models including load-balanced networks.
The statements of stability are established under quite general
assumptions on arrival and departure processes by reduction to the
corresponding models with autonomous service mechanism.

Now we discuss how these results can be extended to the models of
queues and networks allowing batch arrivals and batch departures.
For this purpose, consider the queueing system with batch arrivals
and departures and autonomous service. For this queueing system
let $\mathscr{A}(t)$ denote arrival process and let
$\mathscr{D}(t)$ departure process, both marked point processes.
(All of the processes considered in this section are assumed to
start at zero.) For the sake of simplicity suppose that the marks
of the point process $\mathscr{D}(t)$ all are of the constant size
$c$ ($c$ is a positive integer number), and therefore
$\mathscr{D}(t)=cD(t)$. Then, the queue-length process $Q(t)$ has
the following representation (see \cite{Abramov (2006a)}):
$$
Q(t)~=~\mathscr{A}(t)-\sum_{i=1}^c\int_0^t\mathbf{1}_{\{Q(s-)\geq
i\}}\mbox{d}D(s).
$$
It was shown in \cite{Abramov (2006a)} that by using Skorokhod's
reflection principle we arrive at equation
\begin{equation}\label{6.1}
Q(t){\buildrel d\over
=}[\mathscr{A}(t)-\mathscr{D}(t)]-\inf_{s\leq
t}[\mathscr{A}(s)-\mathscr{D}(s)],
\end{equation}
which is similar to that of the process with ordinary departures.
The assumption, that the marks of departure process are a constant
$c$, is specific and associated with concrete  models considered
in \cite{Abramov (2006a)}. Representation \eqref{6.1} remains in
force in general, when a departure process is an arbitrary marked
point process with mutually independent identically distributed
marks. Representation \eqref{6.1} is easily generalized to the
case of JS-queue models. Specifically, for the queue-length
process in the $j$th server of the model $\wp_m$ we have the
similar equation
\begin{equation}\label{6.2}
\begin{aligned}
Q^{(j)}(t)&{\buildrel d\over
=}[\mathscr{A}^{(j)}(t)+\mathscr{A}^{\prime(j)}(t)-\mathscr{D}^{(j)}(t)]\\
&\ \ \ -\inf_{s\leq
t}[\mathscr{A}^{(j)}(s)+\mathscr{A}^{\prime(j)}(t)-\mathscr{D}^{(j)}(s)],
\end{aligned}
\end{equation}
where $\mathscr{A}^{\prime(j)}(t)$ is the corresponding notation
for an opportunistic traffic to the $j$th server of the JS-queue
model (see ref. \eqref{4.4} for comparison). Thus, the case of
batch arrivals and departures is a direct extension of the case of
ordinary arrivals and departures, and the conditions for stability
are similar.

\section*{Acknowledgements}
The advice of Professors Rafael Hassin, Robert Liptser, Yuri Suhov
and Gideon Weiss helped very much to substantially improve the
presentation. The research was supported by Australian Research
Council, grant \#DP0771338.

\end{document}